\documentclass[12pt]{article}
\usepackage{enumerate,amsfonts,amsmath}

\newtheorem {lemma}{Lemma}
\newtheorem {cora}{Corollary}
\newtheorem {prop}{Proposition}
\newtheorem {thm}{Theorem}

\def\E{{\mathbb{E}}}
\def\V{\mathbb{V}}
\def\P{{\mathbb P}}
\def\R {{\mathbb{R}}}
\def\Z {{\mathbb{Z}}}
\def \D{\mathbb{D}}

\def\limt{\lim_{t\to\infty}}
\def\limt0{\lim_{t\to 0}}

\def\lf{\lfloor}
\def\rf{\rfloor}
\def\card{\rm{{card}}}
\def\|{\,|\,}

\newcommand{\BBox}{\rule{6pt}{6pt}}
\newcommand\Cox{$\hfill \BBox$ \vskip 5mm}

\def\proof{\noindent{\sf \underline{Proof}.\ }}
\def\bn#1\en{\begin{align*}#1\end{align*}}
\def\bnn#1\enn{\begin{align}#1\end{align}}

\title{On a  coloured tree with non i.i.d.\ random labels}
\author{Skevi Michael\footnote{Department of Mathematics, University of Bristol, BS8 1TW, U.K.}\ \  and Stanislav Volkov$^*$\footnote{Corresponding author. E-mail: s.volkov@bristol.ac.uk}}
\begin{document}
\maketitle

\begin{abstract}
We obtain new results for the probabilistic model introduced
in~\cite{MPV} and~\cite{VS} which involves a $d$-ary regular tree.
All vertices are coloured in one of $d$ distinct colours so that
$d$ children of each vertex all have different colours. Fix $d^2$
strictly positive random variables. For any two connected vertices
of the tree assign to the edge between them {\it a label} which
has the same distribution as one of these random variables, such
that the distribution is determined solely by the colours of its
endpoints. {\it A value} of a vertex is defined as a product of
all labels on the path connecting the vertex to the root. We study
how the total number of vertices with value of at least $x$ grows
as $x\downarrow 0$, and apply the results to some other relevant
models.
\end{abstract}

\section{Introduction}\label{sec_intro}
In~\cite{VS} Volkov  showed how the $5x+1$ problem can
be approximated by a probabilistic model involving a binary tree with randomly labeled edges,
with distributions of the random variables assigned to edges being
determined by their directions, these random variables being independent.

Menshikov et al.~\cite{MPV} studied a similar model, where random
variables assigned to edges of the tree were  dependent both on
the type of parent vertex and the type of the child, as described
below. At the same time, the results in~\cite{MPV} did not give
the answers to all the questions answered in~\cite{VS}, and this
is the purpose of the current paper. We want to stress that
answering these questions is not a straightforward application of
the previous results, but requires some new additional arguments.

Let $d\ge 2$. We consider the $d$-ary regular rooted tree $T_d$
with vertex set $\V$  (that is, the tree where every vertex has
degree $d+1$ with the exception of the root, $u_0 \in \V$, which
has degree $d$). For the vertices $u, w \in \V$, the following
quantities are defined:
\begin{itemize}
\item{}$\ell(u)$ is the unique self-avoiding path connecting $u$ to the
root;
\item{} $|u|$ is the number of edges in $\ell(u)$;
\item{}$\V_n=\{u\in\V:\ |u|=n\}$ is  the set of $d^{n}$ vertices that lie at graph-theoretical distance $n$ from the
root;
\item{}$u \sim w$ means that $u$ and $w$ are connected by an edge.
\end{itemize}

Among $d$ distinct colours we arbitrarily choose one to colour the
root. All other vertices are coloured from left to right, so that
all $d$ children of each vertex have different colours. We denote
by $c(u)\in \{1,2,..,d\}$ the colour assigned to the vertex $u$.

Now we assign a random variable (label) to each edge as follows.
First, consider $d^2$ strictly positive and non-degenerate random
variables, $\tilde{\xi}_{ij} $, with $i,j \in \{1,2,...d\}$, of
known joint distribution. Now for $u,w \in \V$ such that $u\sim w$
we assign the random variable,  $\xi_{uw}$ to the undirected edge
$(u,w)\equiv (w,u)$, so that:
\begin{itemize}
\item
for every edge $(u,w)$ such that $u$ is the parent of $w$,
$\xi_{uw} \stackrel{\mathcal{D}}{=}$ $ \tilde{\xi}_{c(u)c(w)}$
where $X\stackrel{\mathcal{D}}{=} Y$ means that $X$ and $Y$ have
the same distribution, and

\item
for any collection of edges of the tree $(u_1, w_1),(u_2,
w_2),...,(u_m, w_m)$, where $u_i$ is the parent of $w_i$ $\forall
i \in {1,2,...m}$ and $u_i \neq u_j$ whenever $i \neq j$, the
random variables $\{\xi_{u_iw_i}\}_{i=1}^m$ are independent.
\end{itemize}

For $u \in \V$, we define value $\xi[u]$ to be the product of all the
random variables assigned to the edges of  $\ell(u)$. The main
object of interest in the present paper is
 \bn
Z(x):= \card \{u \in \V : \xi[u]\ge x \}.
 \en
In~\cite{VS} the ultimate object of interest was the complimentary
quantity $Q(x)= \card \{u \in \V : \xi[u]\le x \}$, however, one
can easily see that these two problems are equivalent once we replace
$\tilde\xi_{ij}$  and $x$ by its inverses $(\tilde\xi_{ij})^{-1}$
and $x^{-1}$ respectively; we have chosen to study $Z$ here in order
to be consistent with notations in~\cite{MPV}.


Similar to~\cite{MPV}, we will randomize the colouring to avoid
the disadvantage of the above colouring method, consisting in the
fact that for different $u, w \in \V_n$ the distribution of
$\xi[u]$ may differ from that of $\xi[w]$. In order to achieve
equality of the distributions of $\xi[u]$ for all $u\in \V_n$, let
the colouring be done recursively for $n=1,2,\dots$ as follows. We
first colour the root in any of the possible $d$ colours; next,
assuming that the vertices up to level $n-1$ (i.e., the vertices
that belong in $\V_1, \V_2,\ldots, \V_{n-1}$) are already
coloured, independently for each $v\in\V_{n-1}$ we colour each of
its children in some colour so that no two children have the same
colour, with all $d!$ colourings of the children of $v$ being
equally likely. As a result, each one of the $(d!)^{d^{n-1}}$
possible colourings of $V_n$ has the same probability.

As before, to each edge $(u, w)$ we assign  a random variable
$\zeta_{uw}$, which distribution satisfies the conditions imposed
on $\xi_{uw}$. Define $\zeta[u]$  in the same way as $\xi[u]$; then it is clear that at
every level $n$ the distribution of the {\it unordered} set
$\{\zeta[u], u \in \V_n\}$ is the same as the distribution of
$\{\xi[u], u \in \V_n\}$. This means that the two models will give
the same results for a number of problems, while the randomized
colouring ensures that for any $u, w \in \V_n$
$\zeta[u]\stackrel{\mathcal{D}}{=} \zeta[w]$, even though
$\zeta[u]$ and $\zeta[w]$ could be dependent. In particular,
$Z(x)= \card \{u \in \V : \zeta[u]\ge x \}.$

\section{Results from~\cite{MPV}}\label{sec_prev}
Let probability $\P $ and expectation $\E$ be with respect to the
measure generated {\it both} by a random colouring
$\textbf{c}=\{c(u), u \in \V\}$ and a random environment
$\zeta=\{\zeta_{uw},\  u,w \in \V\ \text{such that } u\sim w\}$. Define the $d \times d$
matrix $m(s)$, $s \in [0,\infty)$, as
 \bn
m(s) := \left( \begin{array}{cccc}
\E[\tilde{\xi}_{11}]^s & \E[\tilde{\xi}_{12}]^s & \ldots& \E[\tilde{\xi}_{1d}]^s \\
\E[\tilde{\xi}_{21}]^s & \E[\tilde{\xi}_{22}]^s & \ldots& \E[\tilde{\xi}_{2d}]^s \\
\vdots & \vdots & \ddots  & \vdots \\
\E[\tilde{\xi}_{d1}]^s & \E[\tilde{\xi}_{d2}]^s & \ldots&
\E[\tilde{\xi}_{dd}]^s
\end{array} \right).
 \en
Let $\rho (s)$ be its largest eigenvalue, then $\rho(s)$ is
positive by Perron-Frobenius theorem for matrices with strictly
positive entries.

Let $\D = \left\{ s \in \R : \E [ \tilde{\xi}_{ij}]^s <\infty
\quad \forall \quad i,j \in \{1,2,...,d\}\right\}$ and
$\text{Int}(\D)$ be its interior. Assume that the conditions
below are satisfied:
 \bnn\label{conds}
 \begin{array}{rcll}
 [0,1] &\subseteq& \D,&\\
 0 &\in& \text{Int}(\D),&\\
 \E |\log\tilde{\xi}_{ij}| &<& \infty                    &\forall \quad i , j \in\{1,2,...,d\},\\
 \E |\tilde{\xi}_{ij}\log\tilde{\xi}_{ij}| &<& \infty    &\forall \quad i , j \in
 \{1,2,...,d\}.
 \end{array}
 \enn

\begin{thm}[Theorem 2 in~\cite{MPV}]\label{th1}
Suppose $x>0$, 
$$
\lambda=\inf _{s\ge 0} \rho(s)
$$
and conditions (\ref{conds}) are
fulfilled. Then
\begin{enumerate}[(a)]
\item if $\lambda <1$, then $Z(x)< \infty$ a.s.;
\item if $\lambda >1$, then $Z(x)= \infty$ a.s.
\end{enumerate}
\end{thm}
For a vertex $u\in\V_n$, let $u_0,u_1,\dots,u_{n-1},u_n\equiv u$ be the consecutive vertices
of the path $\ell(u)$.
The proof of the above theorem is largely based on the following statement from~\cite{dH}.
\begin{lemma}[Lemma 1 in~\cite{MPV}]\label{lem2} Let $S_n = \sum_{i=1}^n \log (\zeta_{u_{i-1}u_i})$
and $k_n(s) = \left(\E[e^{sS_n}]\right)^{1/n}
=\left(\E\left[\prod_{i=1}^n\zeta_{u_{i-1}u_i}^s\right]\right)^{1/n}$.
Suppose (\ref{conds}) is fulfilled. Then
\begin{enumerate}[(a)]
 \item $ k (s) = \lim_{n \rightarrow \infty} k_n (s) \in [0, \infty]$  exists for  all $s$;
 \item $\Lambda(s) = \log \rho(s) - \log \text{d} = \log k (s) \in (-\infty, + \infty]$ is  convex;
 \item the rate function $\Lambda^* (z) =\sup_{s \ge 0}(sz - \Lambda(s))$, $z \in \R$, is convex, lower semi-continuous
 and differentiable in Int($\D$). Moreover,
 $$
 \Lambda^*(z)=\begin{cases}
  s_0(z)z-\Lambda(s_0(z)),  & \mbox{if }z\ge \Lambda'(0), \\
  0, & \mbox{if }z\le \Lambda'(0),
\end{cases}
$$
where $s_0(z)$ is the solution of equation $z-\Lambda'(s)=0$;
\item for all $a > 0$, $$\lim_{n \rightarrow \infty} \frac{1}{n} \log \P \left( \frac{S_n}{n} \geq \log a \right) = -\Lambda^{*}
(\log a).$$
\end{enumerate}
\end{lemma}

\section{Expectation of $Z(x)$}
Here we will need one additional assumption:
 \bnn\label{conda}
\E \left[\tilde{\xi}_{ij}\right]^s \in \textbf{C}^2(\R_+) \quad
\forall \quad i , j \in \{1,2,...,d\}
 \enn
as functions of $s$, which is required to ensure that ${\Lambda}
\in \textbf{C}^2(\R_+)$. Indeed, the characteristic polynomial
$P(s,\lambda)=\det(m(s)-\lambda I)$ of $m(s)$ can be written as
$$
P(s,\lambda)= \sum _{k=0}^d a_k(s) \lambda ^k .
$$
where $a_k(s)\in \textbf{C}^2(\R_+)$, $k=0,1,...,d$, are its
coefficients and $I$ is $d\times d$ identity matrix.
By the Perron-Frobenius theorem, $\rho(s)$ is a {\it
simple} root of this polynomial, hence it is {\it not} a root of
the polynomial $\frac{\partial P(s,\lambda)}{\partial\lambda}=0$.
Hence
$$
\left.\frac{\partial P(s,\lambda)}{\partial\lambda}\right|_{\lambda=\rho(s)}\ne
0
$$
and by the implicit function theorem  we obtain that $\rho(s)$
is continuously differentiable in $s$ as $a_i(s)$ are, i.e.\
$\rho(s) \in \textbf{C}^2(\R_+)$ and therefore $\Lambda \in
\textbf{C}^2(\R_+)$.

Suppose conditions (\ref{conds}) are fulfilled. By
Theorem~\ref{th1} if $\lambda < 1$ then $Z(x) <\infty$ a.s. Also,
since
 \bn
 \Lambda^{*} (z) &= \sup_{s \ge 0}(sz - \Lambda(s))
  =\sup_{s \ge 0}(sz - \log \rho(s) +\log d)
 \en
we have
 \bn
 \Lambda^{*} (0) &= \sup_{s \ge 0}(- \log {\rho}(s) +\log d)
  = -\log \left(\inf_{s \ge 0} {\rho}(s)\right) +\log d
  =-\log\lambda +\log d.
 \en
Therefore,
 \bnn\label{eqZlambda}
 {\lambda} < 1 \iff {\Lambda}^{*} (0) > \log d.
 \enn

From now on assume that indeed $\lambda<1$ and hence $Z(x)$ is
a.s.\ finite for all $x>0$. Observe that $Z(x)$ increases to
$+\infty$ as $x\downarrow 0$. We are now ready to give the main
theorem describing the asymptotical behaviour of $\E[Z(x)]$, thus
generalizing the result of Theorem~3 in~\cite{VS} to a more
general setup of~\cite{MPV} described above.

\begin{thm}\label{th3}
Suppose that conditions (\ref{conds}) and (\ref{conda}) are
fulfilled, and moreover the following are true:
\begin{itemize}
\item[(A1)] $\lambda<1$;
\item[(A2)] $\mu:=-\Lambda'(0) > 0$ (equivalently, $\rho'(0)<0$).
\end{itemize}
Then
$$
\lim_{t\to\infty} \frac{\log\E\left[Z(e^{-t})\right]}{t} \:
\textrm{ exists and is given by } M=\max _{ u\in [0,\mu]
}\frac{\log d - {\Lambda}^{*}(-u)}{u}.
$$
\end{thm}
\proof Let
$$
f(u) = \frac{\log d -\Lambda^*(-u)}{u}.
$$
By the definition of the rate function ${\Lambda}^{*}(z)\ge 0$ for
all $z \in \R$, and also ${\Lambda}^{*}(-\mu)=0$. Since
${\Lambda}^{*}$ is a differentiable and convex function we have
${\Lambda^*}'(-\mu)\equiv
\left.\frac{d\Lambda^*(z)}{dz}\right|_{z=-\mu}=0$. Also
 \bn
 \lim_{u \to +0}f(u) &= -\infty \quad(\textrm{because of A1 and (\ref{eqZlambda})});\\
  f(\mu) &= \frac{\log d- \Lambda^{*}(-\mu)}{\mu}=\frac{\log d}{\mu}>0 \quad(\textrm{because of A2}); \\
  f'(\mu) &= \frac{\mu \cdot {\Lambda^*}'(-\mu)- \log d +
  \Lambda^*(-\mu)} {\mu^2} = -\frac{\log d} {\mu^2}< 0.
 \en
We conclude that $\max _{x \in [0, -{\Lambda}^{'}(0)]} f(x)$
exists and is achieved strictly inside the interval $ (0,\mu)$.
Let $u^*\in(0,\mu)$ denote the point where the maximum of $f(u)$
is achieved.


Keeping in mind that $\zeta[\cdot]$ is the same for the vertices
which appear at the same level of the tree, we derive an
expression for $\E[Z(e^{-t})]$ similar to~\cite{VS}:
 \bn
 \E[Z(e^{-t})] & =\sum_{u \in \V} \P  \left(\zeta [u] \ge e^{-t}\right)=\sum_{n=0}^{\infty}\sum_{u \in \V_n} \P  (\zeta [u] \ge e^{-t})
  = \sum_{n=0}^{\infty}\sum_{u \in \V_n} \P  \left(\log \zeta [u] \ge -t\right)
  \\
 & =\sum_{n=0}^{\infty} d^n \cdot \P  \left(\log \left(\prod_{i=1}^{n}\zeta_{u_{i-1}u_i}\right)\ge -t\right)
  =  \sum_{n=0}^{\infty} d^n \cdot \P  \left(S_n\ge  -t\right),
 \en
where
$$
 S_n  = \sum_{i=1}^{n}\log \left(\zeta_{u_{i-1}u_i}\right).
$$
Hence
 \bn
 \E [Z(e^{-t})] &=\sum_{n=0}^{\infty} \exp\left\{ n\log d + \log \P  \left({S}_n \ge -t\right) \right\}
 =\sum_{n=0}^{\infty} e^{t U_n }
 \en
where
$$
 U_n = \frac{ \log d +  \frac 1n \log \P\left({{S}_n}/{n}\ge -t/n\right)}{t/n}.
$$
First we get the upper bound for $\E [Z(e^{-t})]$.


By Lemma~\ref{lem2}  $\Lambda^*$ is a continuous function and
$\Lambda^*(0)>\log d$, therefore, there are $\epsilon\in(0,\mu)$
and $\bar\delta>0$ such that for all $\delta\in (0,\bar\delta)$ we
have $\Lambda^*(-\epsilon)>\log d+2\delta$. In turn, by part (d)
of Lemma~\ref{lem2} there is an $n_0=n_0(\epsilon,\delta) \in
\mathbb{N}$ such that for all $n\ge n_0$
 \bnn\label{eqU1}
  \frac{1}{n} \log \P \left({{S}_n}/{n} \ge -\epsilon \right) \le   -\Lambda^* (-\epsilon)+\delta \le -(\log d + \delta).
 \enn
%
%
On the other hand, when $n\ge  t/\epsilon$
 \bnn\label{eqU2}
      \P\left(S_n/n \ge -\epsilon\right)
  \ge \P\left(S_n/n\ge-t/n\right).
 \enn
Plugging the inequalities (\ref{eqU1}) and  (\ref{eqU2}) into the
expression for $U_n$ for $n\ge \max{\{n_0, t/\epsilon\}}$ we
obtain $U_n \leq -\frac{n\delta}{t}$. Assume that $t$ is
sufficiently large. Then $t/\epsilon>n_0$ yielding
 \bnn\label{eq1}
   \sum_{n=\lf \frac{t}{\epsilon} \rf +1}^{\infty} e^{tU_n} \leq \sum_{n=0}^{\infty} e^{-n\delta}
  =  \frac{1}{1-e^{-\delta}}.
 \enn
%
%
Secondly,
 \bnn\label{eq2}
 \sum_{n=0}^{\left\lf \frac{t}{\mu}\right\rf} d^n \cdot \P \left(S_n \ge -t\right)
 & \leq  \sum_{n=0}^{\left\lf \frac{t}{\mu}\right\rf} d^n
 \leq \left(\left\lf \frac{t}{\mu}\right\rf +1 \right)  e^{\frac{t\log d }{\mu}}
   \leq  \left(\frac{t}{\mu} +1 \right)
   e^{tM}
 \enn
since $\frac{\log d}{\mu}=f(\mu)\le M$.

To complete the first part of the proof for the upper bound, we
need to study the case when
 \bnn\label{eqMid}
  n\in \left[\frac{t}{\mu} ,\frac{t}{\epsilon}\right] \iff \frac
  tn\in[\epsilon,\mu].
 \enn
%
%
The proof of the following statement is deferred until
Section~\ref{AppA}.
\begin{prop}\label{prop4}
Let $a_1, a_2 \in \R$  be such that $a_1 < a_2$. Then for any
$\delta>0$ there is an $n_1=n_1(a_1,a_2,\delta)$ such that
$$
\frac{1}{n} \log \P \left(\frac{{S}_n}{n} \geq a\right)\leq
-{\Lambda}^{*}(a)+ \delta \quad \text{for all $a \in [a_1 ,a_2]$
and $n\ge n_1$}.
$$
\end{prop}

Set
$$
 \quad a=-\frac{t}{n},\quad a_1=-\mu, \quad
 a_2=-\epsilon.
$$
Note that (\ref{eqMid}) implies $a\in[a_1,a_2]$, hence the
conditions of Proposition~\ref{prop4} are fulfilled, as long as
$t$ is large enough, namely $ t > \mu n_1$. Consequently,
 $$
 \frac{1}{n} \log \P  \left({S_n}/{n} \geq -t/n\right)  \leq  -{\Lambda}^{*}\left(-t/n\right)+ \delta
 $$
yielding
 $$
 U_n \leq \frac{\log d -{\Lambda}^{*}\left(-t/n\right)+\delta} {t/n}  \leq f(t/n) + \frac{n\delta}{t}
 \leq M + \delta/\epsilon
 $$
since $t/n$ satisfies (\ref{eqMid}). As a result
 \bnn\label{eq3}
 \sum_{n=\left\lf \frac{t}{\mu}\right\rf +1}^{\left\lf \frac{t}{\epsilon}\right\rf} e^{tU_n}
 \leq \frac{t}{\epsilon} \cdot e^{t (M+\delta/\epsilon)}.
 \enn
Consequently, combining (\ref{eq1}), (\ref{eq2}) and (\ref{eq3})
together for $t$ sufficiently large we can obtain the upper bound
as follows:
 \bnn\label{eq4}
 \E[Z(e^{-t})] & =\sum_{n=0}^{\left\lf \frac{t}{\mu}\right\rf} e^{t  U_n}
  + \sum_{n=\left\lf \frac{t}{\mu}\right\rf +1}^{\left\lf\frac{t}{\epsilon}\right\rf} e^{t U_n}
  + \sum_{n=\left\lf \frac{t}{\epsilon} \right\rf +1}^{\infty}e^{t U_n}
 \nonumber\\ \nonumber\\
 & \leq \left(\frac{t}{\mu} +1 \right) e^{tM}
  +\left(\frac{t}{\epsilon} \right) e^{t(M + \delta /\epsilon)}+\frac{1}{1-e^{-\delta}}
 \nonumber\\ \nonumber\\
  &= C(t,\epsilon,\mu,\delta,M) \,\epsilon^{-1}\, t e^{t(M + \delta /\epsilon)}
 \enn
where
 \bn
  \lim_{t\to\infty} C(t,\epsilon,\mu,\delta,M)=1
 \en
for all $\delta>0$. Taking the logarithm of~(\ref{eq4}) we obtain
 \bn
 \limsup_{t\to\infty}\frac{\log \left(\E[Z(e^{-t})] \right)}t&\le
 M+\delta/\epsilon
 \en
Thus by letting $\delta \to 0$ we have
 \bnn\label{eqtabsup}
 \limsup_{t\to \infty}\frac{\log \left(\E[Z(e^{-t})] \right)}{t} \leq M.
 \enn


Now, we obtain the lower bound for $\E[Z(e^{-t})]$.  Recall that
$u^*$ is the value such that $f(u^*)=M$. Fix a small $\delta>0$.
By part (d) of Lemma~\ref{lem2} there is $n_2=n_2(\delta)$ such
that for all $n\ge n_2$
 \bnn\label{eqstar}
\frac{1}{n} \log \P \left( \frac{{S}_n}{n} \ge -u^* \right) \geq
-{\Lambda}^{*} (-u^*)-\delta.
 \enn
For any $t>n_2 u^*$ define $n^*=n^*(t)=\lf t/u^*\rf\ge n_2$. Then
$t/n^*\ge u^*$, moreover $t/n^*=u^*[1+O(1/t)]$. Therefore,
using~(\ref{eqstar}) we obtain
 \bn
 U_{n^*}&\ge  \frac{ \log d +  \frac 1{n^*} \log \P\left(S_{n^*}/n^* \ge -u^*\right)}{t/n^*}
 \ge \frac{ \log d   -\Lambda^*(-u^*)-\delta}{u^*[1+O(1/t)]}\\
 & = M-\delta/u^* +O(1/t).
 \en
Recalling
 $$
\E[Z(e^{-t})]=\sum_{n=0}^{\infty} e^{tU_n}\ge e^{t U_{n*}}
 $$
we obtain
$$
 \liminf_{t\to\infty} \frac{\log \E[Z(e^{-t})]}{t} \ge M-\delta/u^*.
$$
Since $\delta>0$ is arbitrary, this yields
$\liminf_{t\to\infty}\frac{\log \E[Z(e^{-t})]}{t} \ge M $ which,
together with~(\ref{eqtabsup}), concludes the proof. \Cox

In fact, the result of Theorem~\ref{th3} can be rewritten in a
somewhat simpler form.

\begin{cora}\label{corrr}
Suppose that all the assumptions made in Theorem~\ref{th3} hold.
Then
$$
\lim_{t\to\infty}\frac{\log\E\left[Z(e^{-t})\right]}{t}=\min \{s
\in {\D}:{\rho}(s)=1\}.
$$
\end{cora}
Before we present the proof, observe that $\rho(0)=d\ge 2$ and
$\inf_{s\ge 0} \rho(s)  \equiv \lambda<1$, hence $\min \{s \in
{\D}:{\rho}(s)=1\}$ is well defined.

 \proof
Form Lemma~\ref{lem2}, part (b), it follows that we only need to
show that
$$
\min \{s \in {\D}:{\Lambda}(s)=-\log d\}=M
$$
where $M$ is defined in the statement of Theorem~\ref{th3}.

By Lemma~\ref{lem2}, part~(c),
 \bnn\label{eq8a}
 \Lambda^*(z)=z s_0(z)-\Lambda(s_0(z))
 \enn
where $s_0(z)$ solves $ \Lambda'(s_0(z))=z.$ Note that
$s_0(z)=(\Lambda')^{-1}(z)$ is uniquely defined, since $\Lambda$
is strictly convex due to non-degeneracy assumptions
(see~\cite{MPV}, Section~5.4, right after formula~(5.10) there),
yielding that $\Lambda'(s)$ is strictly increasing.
Since $\Lambda'(s)\in{\mathbf C}(\R_+)$ from the arguments after
equation (\ref{conda}), we conclude that  $s_0(z)$ is continuously
differentiable and increasing in~$z$. This implies
 \bnn\label{eq8b}
 {\Lambda^*}'(z)=s_0(z) \quad\text{for all } z.
 \enn

Recall that
$$
f(u) = \frac{\log d -{\Lambda}^{*}(-u)}{u}
$$
and $u^*$ is the point where the maximum of $f$ on the segment
$[0,\mu]$ is achieved; in the proof of Theorem~\ref{th3} we have
shown that $0<u^*<\mu$. Using (\ref{eq8a})  and (\ref{eq8b}) have
 \bnn\label{eq7}
 f'(u)&=\frac{u{\Lambda^*}'(-u) -\log d + \Lambda^*(-u)}{{u}^2}=
 \frac{u s_0(-u) -\log d +[-u s_0(-u)-\Lambda(s_0(-u)) ] }{{u}^2} \nonumber \\
 &=-\frac{\log d +\Lambda(s_0(-u))}{{u}^2}=\frac{s_0(-u)-f(u)}{u}.
 \enn
We know $\Lambda(0)=0$, and from (A1) it follows that $\inf_{s\ge
0} \Lambda(s)<-\log d$, hence from the strict convexity of
$\Lambda$ it follows the set $\{s\ge 0: \ \Lambda(s)=-\log d\}$
contains either $1$ or $2$ points.  Now, if $0<s_1<s_2$ are such
that $\Lambda(s_1)=\Lambda(s_2)=-\log d$, from the convexity it
follows $\Lambda(s)+\log d>0$ for $s<s_1$ and $s>s_2$, while
$\Lambda(s)+\log d<0$ for $s\in (s_1,s_2)$. Suppose
$s_1=s_0(-u_1)$ and $s_2=s_0(-u_2)$, then $u_1>u_2$ (recall that
$s_0(z)$ is increasing), and $f'(u)<0$ for $u<u_2$ and $u>u_1$
while $f'(u)>0$ for $u\in (u_2,u_1)$. This implies that $u^*=u_1$
is the point where the maximum is really achieved.
On the other hand, from (\ref{eq7}) we see that $f'(u)=0$ implies
$f(u)=s_0(-u)$ thus yielding $M=f(u_1)=s_0(-u_1)=s_1$ which
concludes the proof. \Cox

\section{Applications and remaining proof}
The construction studied in this paper relates to many other
probabilistic models; see~\cite{MPV}. These applications include
random walks in random environment, first-passage percolation,
multi-type branching walks among others. Here, we will only focus
on the two of them for which Theorem~\ref{th3} provides additional
information.

\subsection{First-passage percolation}
Consider the coloured tree $T_d$ as constructed in
Section~\ref{sec_intro}. To each edge $(u,w)$, where $u$ is the
parent of $w$ we assign a random variable $\tau_{uw}$ which
denotes the {\it passage time} from vertex $u$ to vertex $w$ and
can be one of the $d^2$ possible types $\tilde\tau_{ij}$,
$i,j=1,\dots,d$; the type is determined by the colours of the
edge's endpoints. We assume for simplicity that all the passage
times are independent. We want to study
$$
R(t)= \card \{ u \in \V : \sum_{(v,w)\in \ell(u)} \tau_{vw}\leq t
\}
$$
that is, the number of vertices of the tree which can be reached
by a particle traveling at unit speed  by time $t$; as in
Section~5.3 of~\cite{MPV}, we allow the passage times to be
negative, indicating a sort of `speeding up' of the motion.
Proposition~3 in~\cite{MPV} provides a criterion for finiteness of
$R(t)$. Using our Theorem~\ref{th3} and Corollary~\ref{corrr} we
obtain a much finer result:

\begin{prop}\label{prop6}
Let $\tilde\xi_{ij}=e^{-\tilde \tau_{ij}}$, $i,j=1,\dots,d$.
Suppose that $m(s)$, $\rho(s)$, $\D$, and $\lambda$ are the same
as in Section~\ref{sec_prev}. If $\lambda<1$ and $\rho'(0)<0$ then
$$
 \lim_{t\to\infty} \frac{\E [R(t)]}{t}
  =\min \{s \in {\D}:{\rho}(s)=1\}.
$$
\end{prop}

\subsection{Multi-type branching random walks on $\R$}
Suppose there are  $d$ different types of particles and $d^2$
positive random variables, $\tau_{ij}$,   $i,j=1,2,\dots,d$, whose
joint distribution is non-degenerate, and define the following
process on $\R$. The process starts at time $n=0$ with one
particle of type $j \in \{1,2,\ldots,d\}$ located at point $0$,
write this as $X_1^{(0)}=0$. At time $n=1$ this particle splits
into $d$ other particles which have different types and take their
position $X^{(1)}_1,X^{(1)}_2,...,X^{(1)}_d$ on the real line. The
distributions of the jumps $X^{(1)}_k-X^{(0)}_m$ are assumed to be
independent for different $k$'s and $m$'s. Now, at time $n=2$ the
first generation particles split into other particles, following
the same rules as the original particle, giving a total of $d^2$
new particles located somewhere on $\R$. If we let this procedure
to continue, at time $n$ we will get exactly $d^n$ particles with
positions $X^{(n)}_1,X^{(n)}_2,...,X^{(n)}_{d^n} \in \R$. Suppose
that the jump from an ancestor to a descendant, say
$X^{(n)}_k-X^{(n-1)}_m$, has the distribution of $\tau_{ij}$
provided the particle at $X^{(n-1)}_m$ is of type $i$ and  the
particle at $X^{(n-1)}_k$ is of type $j$, thus the jump
distribution depends on the types of both the parent and the
offspring. Such a model was considered in~\cite{BJ}
and~\cite{MPV}.

Again, set $\tilde{\xi}_{ij} = e^{-\eta_{ij}}$ and let $\rho(s)$
and $\lambda$ be the same as in Section~\ref{sec_prev}.

\begin{prop}[Proposition 5 in~\cite{MPV}]\label{prop8}
Let $x_0\in \R$ be the unique solution of the equation $\inf_{s\ge
0} e^{sx_0} \rho(s)=1$. Then
$$
 \lim_{n\to\infty}\frac{\min\{X_k^{(n)},\ k=1,2,\dots,d^n\}}{n}=x_0 \quad\text{a.s.}
$$
\end{prop}

Observe that the definition $\tilde{\xi}_{ij}$ above implies 
that $Z(e^{-t})$ corresponds to the number of particles of all generations 
that lie to the left of $t$. Hence, our Theorem~\ref{th3} and Corollary~\ref{corrr} give the following
result about the expected number of visits to $(-\infty,t]$ by
particles of all generations of our branching random walk:
\begin{prop}\label{prop8bis}
Suppose that $\lambda<1$ and $\rho'(0)<0$. Then
$$
 \lim_{t\to\infty}\frac{\log \left( \E \left[
 \sum_{n=1}^{\infty} \card\left\{i\in\{1,2,\dots,d^n\}: \ X_i^{(n)}\le t \right\}\right]\right)}{t}
  =\min \{s \in {\D}:{\rho}(s)=1\}.
$$
\end{prop}

\subsection{Proof of Proposition~\protect\ref{prop4}}\label{AppA}
Firstly, we know that ${\Lambda}^{*}$ is continuous on a compact
set $[a_1,a_2]$ $\iff$  ${\Lambda}^{*}$ is uniformly continuous on
$[a_1, a_2]$ by uniform continuity theorem.

Fix $\delta > 0$. Then we can choose $\tau >0$ small so that, for
$x^{'}, x^{''} \in [a_1, a_2]$
  \bnn\label{eq15}
\left|\Lambda^*(x')- \Lambda^*(x'')\right|
 \le \frac{\delta}{2} \quad \textrm{whenever}\quad |x'- x''|\le \tau.
 \enn
Then we choose an $m \in \Z$ and a sequence of real numbers $x_1,
x_2, \cdots, x_m$ such that,
\begin{center}
 $a_1=x_1 < x_2 < \ldots < x_{m-1} < x_m =a_2$ and\\
$x_{i+1}-x_{i} <\tau \ \forall \ i \in \{1,2,\cdots,m-1\}.$
\end{center}
By Lemma~\ref{lem2}, for each $i\in\{1,2,\cdots,m\}$ there is an
$n_{i}$ such that
 \bnn\label{eq16}
\frac{1}{n}\cdot \log \P  \left(\frac{{S}_n}{n} \geq x_i\right)
\leq - {\Lambda}^{*}(x_i) + \frac{\delta}{2} \quad \forall n\geq
n_i.
 \enn
Define $\tilde{n}:= \max\{n_1, n_2, \cdots, n_m\}<\infty$.

Now, $\forall \: a \in (a_1, a_2)$ there is a
$j\in\{1,2,\cdots,m-1\}$ such that $x_j\le a \le x_{j+1}$.
Consequently, for all $n\ge \tilde n$
 \bn
\frac{1}{n}\cdot \log \P  \left(\frac{{S}_n}{n} \geq a\right)
 &\leq \frac{1}{n}\cdot \log \P  \left(\frac{{S}_n}{n} \geq x_j\right)
  \stackrel{(\text{by }\ref{eq16})}{\leq}
   - \Lambda^*(x_j) +  \frac{\delta}{2}
 \\ & \leq
  \left[ - \Lambda^*(a) + \frac{\delta}{2}\right] +\frac{\delta}{2} = -\Lambda^*(a) +
\delta.
 \en
where the final inequality follows from (\ref{eq15}) and the fact
that $|a-x_{j+1}|<\tau$. \Cox

\end{document}